\documentclass[titlepage,12pt]{article} 

\usepackage{graphicx}
\usepackage{amssymb,amsthm,amsmath} 
\usepackage{mathrsfs}
\usepackage[a4paper]{geometry}
\usepackage{datetime2}
\usepackage[utf8]{inputenc}
\usepackage[italian,english]{babel}
\usepackage{hyperref}

\date{}

\newcommand{\ep}{\varepsilon}
\newcommand{\re}{\mathbb{R}}

\newcommand{\n}{\mathbb{N}}

\newcommand{\rebar}{\overline{\mathbb{R}}}

\newcommand{\Eep}{E_\ep}
\newcommand{\Iep}{\mathcal{I}_\ep}
\newcommand{\Iepn}{\mathcal{I}_{\ep_n}}
\newcommand{\uep}{u_\ep}

\newcommand{\unk}{u_{n_k}}
\newcommand{\Vep}{V_\ep}
\newcommand{\Vepn}{V_{\ep_n}}
\newcommand{\Vepnk}{V_{\ep_{n_k}}}
\newcommand{\X}{\mathbb{X}}
\newcommand{\LHS}{\operatorname{LHS}}
\newcommand{\RHS}{\operatorname{RHS}}

\newtheorem{thm}{Theorem}[section]
\newtheorem{thmbibl}{Theorem}

\newtheorem{rmk}[thm]{Remark}
\newtheorem{prop}[thm]{Proposition}

\newtheorem{defn}[thm]{Definition}

\newtheorem{lemma}{Lemma}
\newtheorem{open}[thm]{Open problem}

\title{On the elliptic regularization of maximal slope curves}

\author{
Massimo Gobbino\vspace{1ex}\\ 
{\normalsize Università degli Studi di Pisa} \\
{\normalsize Dipartimento di Matematica}\\ 
{\normalsize PISA (Italy)}\\  
{\normalsize e-mail: \texttt{massimo.gobbino@unipi.it}}
\and
Filippo Quattrocchi\vspace{1ex}\\ 
{\normalsize Université Paris-Saclay, CNRS, Laboratoire de mathématiques d’Orsay} \\
{\normalsize ParMA, Inria Saclay}\\ 
{\normalsize ORSAY (France)}\\
{\normalsize e-mail: \texttt{filippo.quattrocchi@universite-paris-saclay.fr}}
}

\begin{document}
\maketitle

\begin{abstract}

We consider the so-called Weighted Energy-Dissipation approximation of gradient flows in metric spaces. We prove the convergence of minimizers to maximal slope curves under assumptions that are close to those of the original existence theory for metric gradient flows. The proof relies on an energy-dissipation identity for the approximating curves and on an adaptation of a general stability result for maximal slope curves.
\vspace{6ex}

\noindent{\bf Mathematics Subject Classification 2020 (MSC2020):} 
49J40, 49J45, 35K90.

\vspace{6ex}


\end{abstract}

\section{Introduction}

In this paper we consider the family of functionals
\begin{equation}
    \Iep(u):=\int_0^{+\infty}\left\{\frac{\ep}{2}|u'|(t)^2+F(u(t))\right\}\frac{e^{-t/\ep}}{\ep}\,dt,
    \label{defn:Gep}
\end{equation}
where $\ep>0$, the function $F$ is defined on a suitable metric space~$\X$, and the trajectory $u:[0,+\infty)\to\X$ belongs to a natural class of absolutely continuous curves satisfying the initial condition $u(0)=x_0\in\X$.

We prove that, under rather general assumptions on $\X$, $F$, and $x_0$, the functional $\Iep$ admits at least one minimizer $u_\ep$ for every $\ep>0$. Moreover, the family $\{u_\ep\}$ converges, up to subsequences, to a limit trajectory $u$ as $\ep\to0^+$. We also show that every such limit is a maximal slope curve for $F$, namely a trajectory that, roughly speaking, follows the direction of steepest descent of the functional.

The notion of maximal slope curve, as a metric-space counterpart of gradient flows, was introduced by E.~De~Giorgi and collaborators in~\cite{DG:CMP}, and subsequently developed in a series of papers (see, for example,~\cite{MarinoSacconTosques89}). We refer to~\cite{AGS} for a modern systematic presentation of the theory.

A variational approach closely related to~(\ref{defn:Gep}) was introduced by E.~De Giorgi in~\cite{DG:Duke} as an approximation scheme for certain hyperbolic PDEs. In that context, the functional~(\ref{defn:Gep}) is replaced by a second-order analogue involving the term $\ep^2|u''|(t)^2$ instead of $\ep|u'|(t)^2$. De Giorgi's construction was inspired by the elliptic regularization method developed by T.~Ilmanen in~\cite{1994:ilmanen} in his proof of existence for the Brakke mean curvature flow. The second-order functional was later investigated by P.~Tilli and E.~Serra in~\cite{serra_tilli_2016}.

The first-order functional~(\ref{defn:Gep}) was considered independently by the second author in his bachelor's thesis~\cite{quattrocchi} and by R.~Rossi, G.~Savaré, A.~Segatti and U.~Stefanelli in~\cite{rossi_savare_segatti_stefanelli_2019}. The latter paper appeared during the preparation of the thesis. Although the two works address the same convergence problem, they rely on different assumptions and techniques. 

The present paper develops the approach originally pursued in the thesis, which is closer in spirit to the classical theory of maximal slope curves: the main difference from~\cite{rossi_savare_segatti_stefanelli_2019} is not the variational scheme itself, but the framework in which the limit is identified. In the present work, maximal slope curves are understood in the original De Giorgi sense, with the descending slope~$|\nabla F|$ (see Definition~\ref{defn:slope})---as opposed to its lower-semicontinuous relaxation---and by making direct (semi)continuity assumptions on the energy and its slope (see~(\ref{hp:F-cont}) below), rather than under the requirement that the relaxed slope be an upper gradient, cf.~\cite[Theorem 3.6]{rossi_savare_segatti_stefanelli_2019}. Our strategy is first to establish an energy-dissipation identity for the approximating curves, analogous to the one satisfied by curves of maximal slope. We then use the (semi)continuity assumptions on the energy and its slope to pass to the limit in this identity via an adaptation of the
abstract convergence argument of~\cite[Theorem~3.1]{SlowPM2TV}.

We refer to~\cite{serra_tilli_2016,rossi_savare_segatti_stefanelli_2019,stefanelli_2025} for a more extensive bibliography on the subject, as well as for further motivations and possible future developments. Here we only emphasize that the variational procedure considered in this paper provides an alternative to the classical minimizing movement approach for approximating maximal slope curves.

\paragraph{\textmd{\textit{Structure of the paper}}}

This paper is organized as follows. In Section~\ref{sec:statements} we fix the notation, recall the notion of maximal slope curves in metric spaces, and state our assumptions on the function~$F$ together with the main result. In Section~\ref{sec:existence} we recall a compactness criterion for $\X$-valued trajectories and prove the existence of minimizers for~(\ref{defn:Gep}) among all trajectories with a prescribed initial condition. In Section~\ref{sec:compactness} we prove that any family of minimizers is relatively compact. In Section~\ref{sec:limit} we prove that every limit of minimizers is a maximal slope curve. This requires an analysis of the asymptotic behavior of the so-called value function, defined in~(\ref{defn:vf}). Finally, in Section~\ref{sec:open} we collect some open questions naturally arising from the proof.

%\clearpage

\setcounter{equation}{0}
\section{Notation and main result}\label{sec:statements}

To fix notation, we briefly recall the main ingredients of the classical theory of gradient flows in metric spaces. In the sequel, $\rebar := \re \cup \{\pm\infty\}$ denotes the extended real line.

\begin{defn}[Descending metric slope]\label{defn:slope}
\begin{em}

Let $(\X,d)$ be a metric space, let $F:\X\to\rebar$ be a function, and let $x_0\in\X$.

The \emph{descending slope} of $F$ at $x_0$, denoted by $|\nabla F|(x_0)$, is defined to be $+\infty$ whenever~$F(x_0)\notin\re$, and otherwise by
\begin{equation}
    |\nabla F|(x_0):=
    \limsup_{r\to 0^+}
    \frac{F(x_0)-\inf\{F(x):d(x,x_0)\leq r\}}{r}.
    \nonumber
\end{equation}

\end{em}    
\end{defn}
\begin{rmk}[Equivalent definitions of the metric slope]
\begin{em}

The descending metric slope is also known in the literature as the \emph{local (descending) slope} (see, for example,~\cite{DG:CMP,rossi_savare_segatti_stefanelli_2019}). It is usually defined by
\[
    |\nabla F|(x_0)
    :=
    \limsup_{x\to x_0}
    \frac{\max\{F(x_0)-F(x),0\}}{d(x,x_0)}
\]
whenever $F(x_0)$ is finite and $x_0$ is not an isolated point of $\X$. At isolated points in the domain of~$F$, the slope is defined to be zero, while it is set equal to $+\infty$ when $F(x_0)$ is not finite.

It is straightforward to verify that this definition is equivalent to the one given above in Definition~\ref{defn:slope}, which follows~\cite{antonucci}. We prefer the latter formulation because it automatically covers isolated points and avoids the need for an explicit positive part in the numerator.

\end{em}
\end{rmk}

\begin{defn}[Absolutely continuous functions and metric derivative]
\begin{em}

Let $(\X,d)$ be a metric space, and let $[a,b]\subseteq\re$ be an interval. The space $AC^2([a,b];\X)$ consists of all functions $u:[a,b]\to\X$ for which there exists a function $g\in L^2((a,b))$ such that
\begin{equation}
    d(u(t_2),u(t_1))\leq\int_{t_1}^{t_2}g(t)\,dt
    \qquad
    \forall a\leq t_1\leq t_2\leq b.
    \label{defn:AC2}
\end{equation}

For every such function $u$, there exists a minimal function $g$ for which~(\ref{defn:AC2}) holds; see~\cite[Theorem~1.1.2]{AGS}. This minimal function is called the \emph{metric derivative} of $u$ and is denoted by $|u'|(t)$.

We denote by $AC^2_\mathrm{loc}([0,+\infty);\X)$ the set of all functions $u:[0,+\infty)\to\X$ whose restriction to $[0,T]$ belongs to $AC^2([0,T];\X)$ for every $T>0$.

\end{em}    
\end{defn}

\begin{defn}[Maximal slope curve]\label{defn:MSC}
\begin{em}

Let $(\X,d)$ be a metric space, and let $F:\X\to\rebar$ be a function. A \emph{maximal slope curve} for $F$ is a triple $(u,\psi,N)$ such that
\begin{itemize}
    \item $u\in AC^2_\mathrm{loc}([0,+\infty);\X)$;

    \item $\psi:[0,+\infty)\to\re$ is a nonincreasing function satisfying
    \begin{equation}
        \frac{1}{2}\int_{t_1}^{t_2}|u'|(\tau)^2\,d\tau+
        \frac{1}{2}\int_{t_1}^{t_2}|\nabla F|(u(\tau))^2\,d\tau
        \leq
        \psi(t_1)-\psi(t_2)
    \nonumber
    \end{equation}
    for every pair of real numbers $0\leq t_1\leq t_2$;

    \item $N\subseteq(0,+\infty)$ is a set of Lebesgue measure zero such that
    \begin{equation}
        \psi(t)=F(u(t))
        \qquad
        \forall t\in[0,+\infty)\setminus N.
    \nonumber
    \end{equation}
\end{itemize}

\end{em}    
\end{defn}

\begin{rmk}[H\"older continuity of $AC^2$ functions]
\begin{em}
    We observe that functions in $AC^2([0,T];\X)$ are H\"older continuous of order~$1/2$. Indeed, from~(\ref{defn:AC2}) and H\"older's inequality it follows that
    \begin{equation}
        d(u(t_2),u(t_1))\leq
        \left\{\int_{t_1}^{t_2}|u'|(t)^2\,dt\right\}^{1/2}|t_2-t_1|^{1/2}
        \qquad
        \forall\,0\leq t_1\leq t_2\leq T.
        \label{th:AC-Holder}
    \end{equation}

    We also point out that, in Definition~\ref{defn:MSC}, we require $0\notin N$. As a consequence,
    \begin{equation}
        \psi(0)=F(u(0)).
    \nonumber
    \end{equation}
\end{em}
\end{rmk}

In the following definition we collect the assumptions on the function $F$ that will be used throughout the paper.

\begin{defn}[Assumptions on $F$]
\begin{em}

Let $(\X,d)$ be a metric space, and let $F:\X\to\rebar$ be a function. We consider the following assumptions on $F$.
\begin{enumerate}
\renewcommand{\labelenumi}{(\arabic{enumi})}
    \item (Sign condition)
    \begin{equation}
        F(x)\geq 0
        \qquad
        \forall x\in\X.
        \label{hp:F-sign}
    \end{equation}

    \item (Lower semicontinuity) For every convergent sequence $x_n\to x_\infty$ in $\X$, one has
    \begin{equation}
        \liminf_{n\to +\infty}F(x_n)\geq F(x_\infty).
        \label{hp:F-lsc}
    \end{equation}

    \item (Local coercivity) For every $x_0\in\X$, every $r>0$, and every $M\in\re$, the set
    \begin{equation}
        \{x\in\X: F(x)\leq M,\ d(x,x_0)\leq r\}
        \label{hp:F-loc-coerc}
    \end{equation}
    is relatively compact, namely every sequence in this set admits a subsequence converging in $\X$.

    \item (Continuity properties on bounded energy-slope sequences) For every convergent sequence $x_n\to x_\infty$ in $\X$, the following implication holds:
    \begin{equation}
    \fbox{$\begin{array}{c}
        \displaystyle{\sup_{n \in \n}|F(x_{n})|< +\infty} \\
        \noalign{\vspace{1ex}}
        \displaystyle{\sup_{n \in \n}|\nabla F|(x_{n})< +\infty}
    \end{array}$}
    \quad\Longrightarrow\quad
    \fbox{$\begin{array}{c}
        \displaystyle{\lim_{n\to +\infty}F(x_{n})=F(x_\infty),} \\
        \noalign{\vspace{1ex}}
        \displaystyle{\liminf_{n\to +\infty}|\nabla F|(x_{n}) \geq |\nabla F|(x_\infty).}
    \end{array}$}
    \label{hp:F-cont}
    \end{equation}
   \end{enumerate}

\end{em}    
\end{defn}

The assumptions stated above are close to those appearing in the classical existence theory for maximal slope curves, which we now recall. The following result coincides with~\cite[Theorem~2.3.9]{antonucci}; see also~\cite{DG:CMP,AGS,SlowPM2TV,MarinoSacconTosques89} and the references therein.

\begin{thmbibl}[Existence of maximal slope curves]\label{thmbibl:CMP}
    Let $(\X,d)$ be a metric space, let $F:\X\to\rebar$ be a function satisfying~(\ref{hp:F-sign}) through~(\ref{hp:F-cont}), and let $x_0\in\X$ be such that $F(x_0)<+\infty$.

    Then there exists at least one maximal slope curve $(w,\psi,N)$ for $F$ with initial condition $w(0)=x_0$. In addition, every such triple has the following properties:
    \begin{itemize}

        \item the function $\psi$ satisfies $\psi(0)=F(x_0)$ and
        \begin{equation}
            0\leq\psi(t)\leq F(x_0)
            \qquad
            \forall t\geq 0;
            \nonumber
        \end{equation}
    
        \item the metric derivative of $w$ coincides almost everywhere with the slope, namely
        \begin{equation}
            |w'|(t)=|\nabla F|(w(t))
            \qquad
            \text{for almost every } t\geq 0;
            \nonumber
        \end{equation}

        \item $w$ is H\"older continuous of order~$1/2$, and more precisely
        \begin{equation}
            d(w(t_2),w(t_1))\leq
            [2F(x_0)]^{1/2}\cdot|t_2-t_1|^{1/2}
            \qquad
            \forall\, 0\leq t_1\leq t_2.
            \nonumber
        \end{equation}

    \end{itemize}
\end{thmbibl}

%\clearpage

The main result of this paper is that the elliptic regularization scheme based on the family of functionals~(\ref{defn:Gep}) converges to maximal slope curves under the same assumptions that ensure their existence (with the addition of~$(\X,d)$ being complete; see also Open~problem~\ref{open:compl}).

\begin{thm}[Convergence of elliptic regularizations to maximal slope curves]\label{thm:main}
    Let $(\X,d)$ be a complete metric space, let $F:\X\to\rebar$ be a function, and let $\{\Iep\}$ be the family of functionals defined by (\ref{defn:Gep}). 

    Let us assume that $F$ satisfies the sign, lower semicontinuity, and local coercivity conditions (\ref{hp:F-sign}), (\ref{hp:F-lsc}), and (\ref{hp:F-loc-coerc}). Let $x_0\in\X$ be a point with $F(x_0)<+\infty$.

    Then the following statements hold true.
    \begin{enumerate}
    \renewcommand{\labelenumi}{(\arabic{enumi})}
        
        \item \label{statement:existence} \emph{(Existence of minimizers).}
        For every $\ep>0$, the minimization problem
        \begin{equation}
            \min\left\{\Iep(u):
            u\in AC^2_\mathrm{loc}([0,+\infty);\X),\
            u(0)=x_0\right\}
            \label{defn:min-datun}
        \end{equation}
        admits a solution.
        
        \item \label{statement:compactness} \emph{(Compactness of minimizers).}
        Let~$\{\ep_n\}$ be a bounded sequence of positive real numbers and, for every $n\in\n$, let $u_n$ be a minimizer of~(\ref{defn:min-datun}) with $\ep:=\ep_n$.
        
        Then there exist $u_\infty\in AC^2_\mathrm{loc}([0,+\infty);\X)$ and an increasing sequence $\{n_k\}$ of positive integers such that
        \begin{equation}
            \forall T>0
            \qquad
            u_{n_k}\to u_\infty
            \quad
            \text{uniformly on $[0,T]$}.
            \nonumber
        \end{equation}
        
        \item \label{statement:limit} \emph{(Limits of minimizers are maximal slope curves).}
        Let us assume in addition that $F$ satisfies~(\ref{hp:F-cont}) and that $\ep_n\to 0^+$. Let $u_{n_k}\to u_\infty$ be as in the previous statement.
    
        Then there exist a function $\psi_\infty$ and a set $N_\infty$ such that $(u_\infty,\psi_\infty,N_\infty)$ is a maximal slope curve for~$F$ with initial datum $x_0$.
        
    \end{enumerate}
\end{thm}

\begin{rmk}[On the sign condition on $F$]
\begin{em}
The sign condition~(\ref{hp:F-sign}) is merely a normalization. Indeed, adding a constant to $F$ translates the functional $\Iep$ by the same constant, without affecting its minimizers. Moreover, the descending slope remains unchanged, and hence the maximal slope curves are preserved, up to the corresponding translation of the energy representative. Therefore, all our results remain valid for energies bounded from below.

More generally, all results remain true, up to restricting to small values of~$\ep$, under the assumption that $F$ is bounded from below by a negative quadratic function, namely
\begin{equation}
    F(x)\geq -A-Bd(x,z_0)^2
    \qquad
    \forall x\in\X,
    \label{hp:F>=quad}
\end{equation}
for some nonnegative real numbers $A$ and $B$, and some $z_0\in\X$. We refrain from working in this generality, as it only complicates the proofs without introducing new ideas.

On the other hand, when~(\ref{hp:F>=quad}) is not satisfied, there are simple examples in which a maximal slope curve does not exist, at least globally in time.
\end{em}
\end{rmk}

\begin{rmk}[On the initial datum]
\begin{em}

Again for the sake of brevity, we assume that the initial datum is the same point $x_0\in\X$ for every $\ep>0$. The result can be extended in a straightforward way to families of initial data $x_\ep\to x_0$ such that $F(x_\ep)\to F(x_0)<+\infty$.
    
\end{em}
\end{rmk}

\setcounter{equation}{0}
\section{Existence of minimizers}\label{sec:existence}

The proof of Statement~(\ref{statement:existence}) of Theorem~\ref{thm:main} is a straightforward application of the direct method in the calculus of variations, once the following compactness result has been established; see also~\cite[Theorem~2.6]{rossi_savare_segatti_stefanelli_2019} for a related result in a more general metric-topological framework. 
We emphasize that the following compactness criterion assumes $F$ to be nonnegative, measurable, and locally coercive, but does not require lower semicontinuity.

\begin{prop}[Compactness criterion for $AC^2$ trajectories]\label{prop:cpt}
    Let $(\X,d)$ be a complete metric space, let $F:\X\to\rebar$ be a function, and let $\{u_n\}\subseteq AC^2_\mathrm{loc}([0,+\infty);\X)$.

    Let us assume that
    \begin{enumerate}
    \renewcommand{\labelenumi}{(\roman{enumi})}
        \item $F$ is measurable, and it satisfies the sign condition~(\ref{hp:F-sign}) and the local coercivity condition~(\ref{hp:F-loc-coerc});

        \item there exists $x_0\in\X$ such that $u_n(0)=x_0$ for every $n\in\n$;

        \item for every $T>0$, there exists a constant $M_T$ such that
        \begin{equation}
            \int_0^{T}\left\{|u_n'|(t)^2+F(u_n(t))\right\}\,dt\leq M_T
            \qquad
            \forall n\in\n.
            \label{hp:CPT-bound}
        \end{equation}
    \end{enumerate}
       
    Then the sequence~$\{u_n\}$ is relatively compact with respect to uniform convergence on compact subintervals of $[0,+\infty)$, namely there exist $u_\infty\in AC^2_\mathrm{loc}([0,+\infty);\X)$ and an increasing sequence $\{n_k\}$ of positive integers such that
    \begin{equation}
        u_{n_k}\to u_\infty
        \qquad
        \text{uniformly on $[0,T]$ for every $T>0$}.
    \nonumber
    \end{equation}

\end{prop}

\begin{proof}

The proof is a straightforward variant of the classical proof of the Arzelà--Ascoli theorem, so we limit ourselves to outlining the main steps.

Let us fix $T>0$.
\begin{itemize}
    \item Step~1. The family $\{u_n\}$ is equi-H\"older continuous of order~$1/2$ on $[0,T]$. This follows from~(\ref{th:AC-Holder}) and the $L^2$ bound on the metric derivatives provided by~(\ref{hp:CPT-bound}), together with the sign condition~(\ref{hp:F-sign}).

    \item Step~2. There exist a countable dense set $\mathcal{D}_T\subseteq[0,T]$, an increasing sequence $\{n_k\}$ of positive integers, and a function $u_\infty:\mathcal{D}_T\to\X$ such that
    \begin{equation}
        \lim_{k\to +\infty}u_{n_k}(t)=u_\infty(t)
        \qquad
        \forall t\in\mathcal{D}_T.
    \nonumber
    \end{equation}
    This is the crucial step where the local coercivity condition~(\ref{hp:F-loc-coerc}) and the bound on the integrals of $F(u_n(t))$ come into play. We refer to~\cite{rossi_savare_segatti_stefanelli_2019} for the details.

    \item Step~3. The function $u_\infty$ defined in Step~2 is H\"older continuous of order~$1/2$ on $\mathcal{D}_T$, and therefore extends uniquely to a H\"older continuous function on the whole interval $[0,T]$. This is the only point in the proof where the completeness of~$\X$ is essential.

    \item Step~4. One proves that the convergence in Step~2 actually holds for every $t\in[0,T]$, with the same subsequence $\{n_k\}$.

    \item Step~5. One proves that the convergence is in fact uniform on $[0,T]$.
\end{itemize}

Since $T>0$ is arbitrary, a standard diagonal argument yields a single subsequence converging uniformly on every bounded interval.
\end{proof}    

%\clearpage

We are now ready to prove Statement~(\ref{statement:existence}) of Theorem~\ref{thm:main}. To this end, let $\{u_n\}$ be a minimizing sequence. Since the constant function $u(t)\equiv x_0$ is a competitor for~(\ref{defn:min-datun}) with energy $\Iep(u)=F(x_0)$, we can assume that $\Iep(u_n)\leq F(x_0)$ for every $n\in\n$. Therefore, by the definition of $\Iep$,
\begin{equation}
    \int_0^T |u_n'|(t)^2\,dt\leq
    2\Iep(u_n)e^{T/\ep}\leq
    2F(x_0)e^{T/\ep}
    \nonumber
\end{equation}
and
\begin{equation}
    \int_0^T F(u_n(t))\,dt\leq
    \Iep(u_n)\ep e^{T/\ep}\leq
    F(x_0)\ep e^{T/\ep}
    \nonumber
\end{equation}
for every $T>0$ and every $n\in\n$. Thus the sequence $\{u_n\}$ satisfies the assumptions of Proposition~\ref{prop:cpt}. Hence, up to extracting a subsequence (not relabeled), $u_n$ converges uniformly on every bounded interval $[0,T]$ to some $u_\infty$, and in particular $u_n(t)\to u_\infty(t)$ for every $t\geq0$.

Now we pass to the limit in the two terms in the definition of $\Iep(u_n)$. As for the second term, since $F$ is lower semicontinuous and bounded from below, Fatou's lemma yields
\begin{align}
    \liminf_{n\to +\infty}\int_0^{+\infty}F(u_n(t))\frac{e^{-t/\ep}}{\ep}\,dt
    &\geq
    \int_0^{+\infty}\liminf_{n\to +\infty}F(u_n(t))\frac{e^{-t/\ep}}{\ep}\,dt 
    \notag\\
    &\geq
    \int_0^{+\infty}F(u_{\infty}(t))\frac{e^{-t/\ep}}{\ep}\,dt.
    \label{est:lsc-Fu}
\end{align}

As for the term involving the metric derivative, we claim that
\begin{equation}
    \liminf_{n\to +\infty}\int_0^{+\infty}\frac{1}{2}|u_n'|(t)^2 e^{-t/\ep}\,dt\geq
    \int_0^{+\infty}\frac{1}{2}|u_{\infty}'|(t)^2 e^{-t/\ep}\,dt.
    \label{est:lsc-u'}
\end{equation}

Indeed, let us consider a subsequence $\{n_k\}$ that realizes the liminf. For every $T>0$, the sequence $\{|u_{n_k}'|\}$ is bounded in $L^2((0,T))$, and therefore, up to passing to a further subsequence, not relabeled, it converges weakly to some function $g$ in $L^2((0,T))$.

It is not difficult to verify that
\begin{equation}
    g(t)\geq |u_\infty'|(t)
    \qquad
    \text{for almost every $t\in(0,T)$}
    \nonumber
\end{equation}
and, since multiplication by $e^{-t/(2\ep)}$ defines a bounded linear operator on $L^2((0,T))$,
\begin{equation}
    |u_{n_k}'|(t)e^{-t/(2\ep)}
    \rightharpoonup
    g(t)e^{-t/(2\ep)}
    \qquad
    \text{weakly in $L^2((0,T))$.}
    \nonumber
\end{equation}

Hence, by the lower semicontinuity of the norm,
\begin{align*}
    \liminf_{n\to +\infty}
    \int_0^{+\infty}|u_n'|(t)^2e^{-t/\ep}\,dt
    &=
    \lim_{k\to +\infty}
    \int_0^{+\infty}|u_{n_k}'|(t)^2e^{-t/\ep}\,dt
    \\
    &\geq
    \liminf_{k\to +\infty}
    \int_0^T|u_{n_k}'|(t)^2e^{-t/\ep}\,dt
    \\
    &\geq
    \int_0^Tg(t)^2e^{-t/\ep}\,dt
    \notag\\
    &\geq
    \int_0^T|u_\infty'|(t)^2e^{-t/\ep}\,dt,
\end{align*}
and the conclusion follows from the arbitrariness of $T$.

Finally, from~(\ref{est:lsc-Fu}) and~(\ref{est:lsc-u'}) we deduce that
\begin{equation}
    \liminf_{n\to +\infty}\Iep(u_n)\geq
    \Iep(u_\infty),
    \nonumber
\end{equation}
which shows that $u_\infty$ is a minimizer.
\qed

\setcounter{equation}{0}
\section{Compactness of minimizers}\label{sec:compactness}

In the following result we state two important properties of minimizers of~(\ref{defn:min-datun}), both of which follow from the Euler--Lagrange equation associated with the minimization problem. We point out that the proof of Statement~(\ref{statement:compactness}) of Theorem~\ref{thm:main} is a straightforward consequence of the uniform energy bound~(\ref{th:min-bounded}) below and Proposition~\ref{prop:cpt}.

In the statement we need the so-called \emph{value function}, namely the function $\Vep:\X\to\rebar$ defined, for every $\ep>0$ and every $x\in\X$, by
\begin{equation}
    \Vep(x):=
    \inf\left\{\Iep(u):u\in AC^2_\mathrm{loc}([0,+\infty);\X),\ u(0)=x\right\}.
    \label{defn:vf}
\end{equation}

The sign condition~(\ref{hp:F-sign}) implies that $\Vep(x)\geq0$ for every $x\in\X$.

\begin{prop}[Properties of minimizers]
    Let $(\X,d)$, $F$, $\Iep$, and $x_0$ be as in Statement~(\ref{statement:existence}) of Theorem~\ref{thm:main}.

    Let $\ep>0$, and let $\uep\in AC^2_\mathrm{loc}([0,+\infty);\X)$ be a minimizer of~(\ref{defn:min-datun}).

    Then the function $\uep$ has the following properties.
    \begin{enumerate}
    \renewcommand{\labelenumi}{(\arabic{enumi})}
        \item \emph{(Uniform energy bound)} For every $T>0$ one has
        \begin{equation}
            \int_0^T\left\{|\uep'|(t)^2+F(\uep(t))\right\}dt\leq
            \left(T+1+\frac{\ep}{2}\right)F(x_0).
            \label{th:min-bounded}
        \end{equation}

        \item \emph{(Energy identity of maximal slope curve type)} For every $0\leq t_1<t_2$ one has
        \begin{multline}
            \qquad
            \frac{1}{2}\int_{t_1}^{t_2}|\uep'|(t)^2\,dt+
            \frac{1}{2}\int_{t_1}^{t_2}\frac{2}{\ep}\left[F(\uep(t))-\Vep(\uep(t))\strut\right]dt
            \\
            =\Vep(\uep(t_1))-\Vep(\uep(t_2)).
            \qquad
            \label{th:min-CMP}
        \end{multline}    
    \end{enumerate}
\end{prop}

\begin{proof}

We only sketch the main ideas of the proof, referring to Sections~3 and~4 of~\cite{rossi_savare_segatti_stefanelli_2019} for analogous arguments in a more general setting. Given a minimizer $u_\ep$ of~(\ref{defn:min-datun}), the key idea is to consider the function defined by
\begin{equation}
    \Eep(t):=F(\uep(t))-\frac{\ep}{2}|\uep'|(t)^2
    \qquad
    \forall t\geq 0.
    \nonumber
\end{equation}

\begin{itemize}
    \item Step~1. The function $\Eep$ is locally absolutely continuous, and
    \begin{equation}
        \Eep'(t)=-|\uep'|(t)^2
        \qquad
        \text{for almost every $t\geq 0$.}
    \nonumber
    \end{equation}
    This follows by computing the \emph{inner variation} of the functional $\Iep$ at $\uep$.

    \item Step~2. From Step~1 one immediately deduces that $\Eep$---more precisely, its continuous representative---satisfies
    \begin{equation}
        \Eep(t_1)-\Eep(t_2)
        =
        \frac{1}{2}\int_{t_1}^{t_2}|\uep'|(t)^2\,dt+
        \frac{1}{2}\int_{t_1}^{t_2}\frac{2}{\ep}\left[F(\uep(t))-\Eep(t)\right]dt
        \label{eqn:cmp-Eep}
    \end{equation}
    for every $0\leq t_1\leq t_2$.

    \item Step~3. Since
    \begin{equation}
        |\Eep(t)|\leq F(\uep(t))+\frac{\ep}{2}|\uep'|(t)^2,
    \nonumber
    \end{equation}
    from the definition of $\Iep$ we deduce that
    \begin{equation}
        \int_0^{+\infty}|\Eep(t)|\frac{e^{-t/\ep}}{\ep}\,dt\leq
        \Iep(\uep)<+\infty,
    \nonumber
    \end{equation}
    and therefore
    \begin{equation}
        \liminf_{t\to +\infty}|\Eep(t)|e^{-t/\ep}=0.
    \nonumber
    \end{equation}

    \item Step~4. Again from Step~1 one obtains
    \begin{equation}
        \left(\Eep(t)e^{-t/\ep}\right)'=
        -\left\{\frac{\ep}{2}|\uep'|(t)^2+F(\uep(t))\right\}\frac{e^{-t/\ep}}{\ep}
    \nonumber
    \end{equation}
    for almost every $t\geq 0$, and consequently
    \begin{equation}
        \Eep(T)e^{-T/\ep}-\Eep(S)e^{-S/\ep}=
        \int_T^S\left\{\frac{\ep}{2}|\uep'|(t)^2+F(\uep(t))\right\}\frac{e^{-t/\ep}}{\ep}\,dt
    \nonumber
    \end{equation}
    for every $0\leq T<S$. Thanks to Step~3, letting $S\to+\infty$ along a suitable sequence gives
    \begin{equation}
        \Eep(T)= e^{T/\ep}
        \int_T^{+\infty}\left\{\frac{\ep}{2}|\uep'|(t)^2+F(\uep(t))\right\}\frac{e^{-t/\ep}}{\ep}\,dt.
        \label{eqn:Eep-representation}
    \end{equation}

    \item Step~5. We know that $\Eep$ is nonincreasing by Step~1 and nonnegative by Step~4. Since
    \[
        F(\uep(t))=\Eep(t)+\frac{\ep}{2}|\uep'|(t)^2
        =
        \Eep(t)-\frac{\ep}{2}\Eep'(t),
    \]
    we obtain
    \begin{align*}
        \int_0^T\left\{|\uep'|(t)^2+F(\uep(t))\right\}\,dt
        &=
        \int_0^T\left\{\Eep(t)-\left(1+\frac{\ep}{2}\right)\Eep'(t)\right\}\,dt
        \notag\\
        &\leq
        T\Eep(0)+\left(1+\frac{\ep}{2}\right)\Eep(0)
        \notag\\
        &=
        \left(T+1+\frac{\ep}{2}\right)\Eep(0).
    \end{align*}
    Since $\Eep(0) = \Iep(u_\ep) \le F(x_0)$ by~(\ref{eqn:Eep-representation}) and the minimality of~$u_\ep$, we conclude~(\ref{th:min-bounded}).

    \item Step~6. One proves that $\Eep(t)=\Vep(\uep(t))$ for every $t\geq 0$, so that~(\ref{th:min-CMP}) follows from~(\ref{eqn:cmp-Eep}). This identity follows again from the representation~(\ref{eqn:Eep-representation}). Indeed,~(\ref{eqn:Eep-representation}) can be rewritten as $\Eep(T)=\Iep(w_\ep)$, where $w_\ep:[0,+\infty)\to\X$ is the time translation of $\uep$ defined by $w_\ep(t):=\uep(T+t)$. At this point, the crucial observation is that, if~$\uep$ minimizes $\Iep$ with initial datum $\uep(0)=x_0$, then $w_\ep$ minimizes $\Iep$ with initial datum $\uep(T)$ and therefore 
    $$\Eep(T)=\Iep(w_\ep)=\Vep(\uep(T)).$$ 
    
\end{itemize}

The identity $\Iep(w_\ep)=\Vep(\uep(T))$ is usually referred to as the dynamic programming principle; see~\cite{rossi_savare_segatti_stefanelli_2019}.
\end{proof}

\setcounter{equation}{0}
\section{Passing to the limit}\label{sec:limit}

In order to prove Statement~(\ref{statement:limit}) of Theorem~\ref{thm:main}, the key idea is to pass to the limit in the family of identities~(\ref{th:min-CMP}). These identities have the same structure as the inequalities defining maximal slope curves (cf.~Definition~\ref{defn:MSC}), with the value function $\Vep$ playing the role of the energy, and the term $2[F(\uep(t))-\Vep(\uep(t))]/\ep$ playing the role of the squared slope. This suggests relying on the general principle that ``the limit of gradient flows is the gradient flow of the limit''. A rigorous formulation of this principle was established in~\cite[Theorem~3.1]{SlowPM2TV}, whose statement we recall for the reader's convenience.

\begin{thmbibl}[Passing maximal slope curves to the limit]
    Let $(\X,d)$ be a metric space, let $F:\X\to\rebar$ be a function, and let $F_n: \X \to \rebar$ be a sequence of functions.

    Let us assume that for every $x\in \X$, and every pair of sequences
	$\{n_{k}\}\subseteq\n$ and $\{x_{k}\}\subseteq \X$, we have the
	implication
	\begin{equation}
		\fbox{$\begin{array}{r}
			n_{k}\to +\infty,\ x_{k}\to x  \\
			\noalign{\vspace{1ex}}
			\displaystyle{\sup_{k \in \n} 
			\left\{|F_{n_{k}}(x_{k})|+|\nabla F_{n_{k}}
			|(x_{k})\strut\right\} < +\infty}
		\end{array}$}
		\Longrightarrow
		\fbox{$\begin{array}{l}
			\displaystyle{\lim_{k\to +\infty}F_{n_{k}}(x_{k})=F(x).} \\
			\noalign{\vspace{1ex}}
			\displaystyle{\liminf_{k\to +\infty}|\nabla
			F_{n_{k}}|(x_{k}) \geq |\nabla F|(x).}
		\end{array}$}
		\label{hp:slope-energy}
	\end{equation}
	
	Let $T>0$, and for every $n\in \n$ let $(u_n,\psi_{n},N_{n})$ be
	a maximal slope curve for $F_n$ in $[0,T]$.  Let us assume that
	\begin{equation}
		\sup_{n\in\n}\;\max_{t\in[0,T]}|\psi_{n}(t)|
		< +\infty,
		\label{hp:sup-fn}
	\end{equation}
	and that $\{u_{n}\}$ has a pointwise limit, namely there
	exists $u:[0,T]\to \X$ such that
	$$\lim_{n\to +\infty}u_{n}(t)=u(t)
	\quad\quad
	\forall t\in [0,T] .$$
	
	Then there exist $\psi$ and $N$ such that $(u,\psi,N)$ is a
	maximal slope curve for $F$ (in the sense of~\cite{SlowPM2TV}, see Remark~\ref{rmk:mscin0}).
    \label{thm:lgf-gfl}
\end{thmbibl}

\begin{rmk}[On the initial energy condition]\label{rmk:mscin0}
\begin{em}

    We point out that the notion of maximal slope curve used in~\cite{SlowPM2TV} is slightly weaker than ours. Indeed, in that paper the exceptional set $N$ is allowed to contain the initial time $t=0$, whereas in Definition~\ref{defn:MSC} we require $\psi(0)=F(u(0))$. Consequently, the conclusion of Theorem~\ref{thm:lgf-gfl} does not necessarily guarantee a maximal slope curve in the sense of our definition. However, this additional property is recovered if the approximating curves satisfy our definition and their initial energies converge, namely $F_n(u_n(0))\to F(u(0))$. Indeed, the energy representative~$\psi$ is found in~\cite{SlowPM2TV} as the pointwise limit of a subsequence $\{\psi_{n_k}\}$; hence we obtain
    \[
        \psi(0)=\lim_{k\to+\infty}\psi_{n_k}(0)=\lim_{k\to+\infty}F_{n_k}(u_{n_k}(0))=F(u(0)).
    \]
        
\end{em}
\end{rmk}

In our setting, the quantity $2(F-\Vep)/\ep$ is not known to coincide with the squared metric slope of $\Vep$. Therefore, we cannot apply the preceding theorem directly, but we can adapt its proof strategy to the energy identities~(\ref{th:min-CMP}).

\subsection{Convergence properties of the value function}

To begin with, we need a preliminary result based on maximal slope curves. It relates~$F$, $\Iep$, and $\Vep$, and connects $|\nabla F|$ with the function in the second integral in~(\ref{th:min-CMP}).

\begin{lemma}[Value function and maximal slope curves]\label{lemma:CMP1}
    Let $(\X,d)$ be a metric space, let $F:\X\to\rebar$ be a function that satisfies (\ref{hp:F-sign}) through (\ref{hp:F-cont}), let $x\in\X$ be such that $F(x)<+\infty$, and let $(w,\psi,N)$ be a maximal slope curve for $F$ with initial condition $w(0)=x$.

    Then, for every $\ep>0$, the value function $\Vep$ defined by~(\ref{defn:vf}) satisfies
    \begin{equation}
        \Vep(x)\leq
        \int_0^{+\infty}\left\{\frac{\ep}{2}|w'|(\ep\sigma)^2+\psi(\ep\sigma)\right\}e^{-\sigma}\,d\sigma,
        \label{th:vf-vs-CMP}
    \end{equation}
    and
    \begin{equation}
        \int_0^{+\infty}|\nabla F|(w(\ep\sigma))^2e^{-\sigma}\,d\sigma\leq
        \frac{2}{\ep}\left[F(x)-\Vep(x)\right].
        \label{th:vf-vs-slope}
    \end{equation}
\end{lemma}

\begin{proof}

We observe that the maximal slope curve $w$ is a competitor in the definition of $\Vep(x)$, and hence
\begin{equation}
    \Vep(x)\leq
    \int_0^{+\infty}\left\{\frac{\ep}{2}|w'|(t)^2+F(w(t))\right\}\frac{e^{-t/\ep}}{\ep}\,dt.
    \nonumber
\end{equation}

From the definition of maximal slope curve we know  that $F(w(t))=\psi(t)$ for almost every $t\geq 0$, and at this point the change of variable $t=\ep\sigma$ gives exactly~(\ref{th:vf-vs-CMP}).

In order to prove (\ref{th:vf-vs-slope}), we start by considering the function
\begin{equation}
    h(t):=\frac{1}{2}\int_0^t|w'|(\tau)^2\,d\tau+
    \frac{1}{2}\int_0^t|\nabla F|(w(\tau))^2\,d\tau
    \qquad
    \forall t\geq 0,
    \nonumber
\end{equation}
namely the absolutely continuous function such that $h(0)=0$ and
\begin{equation}
    h'(t):=\frac{1}{2}|w'|(t)^2+\frac{1}{2}|\nabla F|(w(t))^2
    \qquad
    \text{for almost every $t\geq 0$.}
    \label{eqn:h'}
\end{equation}

From the definition of maximal slope curve, we know that
\begin{equation}
    h(t)\leq 
    \psi(0)-\psi(t)=
    F(x)-\psi(t)
    \qquad
    \forall t\geq 0.
    \label{est:h-CMP}
\end{equation}

Since $\psi$ is nonnegative, this implies in particular that $h$ is bounded.  For every $T>0$, an integration by parts gives
\begin{equation}
    \int_0^T h'(\ep\sigma)e^{-\sigma}\,d\sigma=
    \frac{1}{\ep}h(\ep T)e^{-T}
    +\frac{1}{\ep}\int_0^T h(\ep\sigma)e^{-\sigma}\,d\sigma.
    \nonumber
\end{equation}

Since $h$ is bounded, letting $T\to +\infty$ we obtain that
\begin{equation}
    \int_0^{+\infty} h'(\ep\sigma)e^{-\sigma}\,d\sigma=
    \frac{1}{\ep}\int_0^{+\infty} h(\ep\sigma)e^{-\sigma}\,d\sigma.
    \label{est:int-h'}
\end{equation}

Now from (\ref{est:h-CMP}) we know that $h(\ep\sigma)\leq F(x)-\psi(\ep\sigma)$ for every $\sigma\geq 0$, and therefore (\ref{est:int-h'}) implies that
\begin{eqnarray*}
    \int_0^{+\infty} h'(\ep\sigma)e^{-\sigma}\,d\sigma & \leq &
    \frac{1}{\ep}\int_0^{+\infty} \left[F(x)-\psi(\ep\sigma)\right]e^{-\sigma}\,d\sigma
    \\
    & = &
    \frac{1}{\ep}F(x)-\frac{1}{\ep}\int_0^{+\infty}\psi(\ep\sigma)e^{-\sigma}\,d\sigma.
    \nonumber
\end{eqnarray*}

At this point, from (\ref{eqn:h'}) we conclude that
\begin{eqnarray*}
    \int_0^{+\infty}|\nabla F|(w(\ep\sigma))^2e^{-\sigma}\,d\sigma & =  &
    \int_0^{+\infty}\left[2h'(\ep\sigma)-|w'|(\ep\sigma)^2\right]e^{-\sigma}\,d\sigma
    \\
    & \leq &
    \frac{2}{\ep}F(x)-\frac{2}{\ep}\int_0^{+\infty}
    \left\{\frac{\ep}{2}|w'|(\ep\sigma)^2+\psi(\ep\sigma)\right\}e^{-\sigma}\,d\sigma,
\end{eqnarray*}
which proves  (\ref{th:vf-vs-slope}) in view of (\ref{th:vf-vs-CMP}).    
\end{proof}

The following result is the technical core of our argument. Besides the general $\Gamma$-liminf inequality, it provides the convergence of the value functions~(\ref{th:vf-cont}) and the lower semicontinuity estimate for the corresponding pseudo-slopes~(\ref{th:liminf-slope}) under the boundedness assumptions~(\ref{hp:bound-slope}). This statement, a consequence of the continuity and lower semicontinuity properties of~$F$ in~\eqref{hp:F-cont}, allows us to adapt the convergence argument of Theorem~\ref{thm:lgf-gfl} to our setting.

\begin{prop}[Convergence of the value function to $F$]\label{lemma:lim-vf}
    Let $(\X,d)$ be a complete metric space, and let $F:\X\to\rebar$ be a function. For every $\ep>0$, let $\Iep$ be the functional defined in~(\ref{defn:Gep}), and let $\Vep:\X\to[0,+\infty]$ be the value function defined in~(\ref{defn:vf}).

    Let $\ep_n\to 0^+$ be a sequence of positive real numbers, and let $x_n\to x_\infty$ be a convergent sequence in $\X$.

    Then the following statements hold.
    \begin{enumerate}
    \renewcommand{\labelenumi}{(\arabic{enumi})}
        \item \label{statement:gammaliminf} \emph{(Gamma--liminf inequality)} If $F$ satisfies the sign, lower semicontinuity, and local coercivity conditions~(\ref{hp:F-sign}),~(\ref{hp:F-lsc}), and~(\ref{hp:F-loc-coerc}), then
        \begin{equation}
            \liminf_{n\to +\infty}\Vepn(x_n)\geq F(x_\infty).
            \label{th:liminf-vf}
        \end{equation}

        \item \label{statement:deluxe} \emph{(Deluxe behavior on sequences with bounded energy and ``slope'')} If, in addition, $F$ satisfies the continuity and lower semicontinuity properties in~(\ref{hp:F-cont}), and the sequences $\{\ep_n\}$ and $\{x_n\}$ satisfy
        \begin{equation}
            \sup_{n\in\n}F(x_n)<+\infty
            \qquad\text{and}\qquad
            \sup_{n\in\n}\frac{2}{\ep_n}\left[F(x_n)-\Vepn(x_n)\right]<+\infty,
            \label{hp:bound-slope}
        \end{equation}
        then
        \begin{equation}
            \lim_{n\to +\infty}\Vepn(x_n)= F(x_\infty)
            \label{th:vf-cont}
        \end{equation}
        and
        \begin{equation}
            \liminf_{n\to +\infty}\frac{2}{\ep_n}\left[F(x_n)-\Vepn(x_n)\right]\geq
            |\nabla F|(x_\infty)^2.
        \label{th:liminf-slope}
        \end{equation}
    \end{enumerate}
\end{prop}

\begin{proof}

The proof is divided into several steps. We first prove the Gamma--liminf inequality for the value functions. Then we introduce suitable maximal slope curves starting from the points $x_n$, and use them to prove the convergence of $\Vepn(x_n)$ to $F(x_\infty)$. Finally, we use the same argument to obtain the Gamma--liminf inequality for the corresponding slope terms.

\paragraph{\textmd{\textit{Gamma--liminf inequality}}}

Since~(\ref{th:liminf-vf}) is trivial when the left-hand side is $+\infty$, we assume, without loss of generality, that the sequence $\{\Vepn(x_n)\}$ is bounded. In this case, for every $n\in\n$, the infimum in the definition of $\Vepn(x_n)$ is attained by the same direct-method argument of Statement~(\ref{statement:existence}) in Theorem~\ref{thm:main}, namely there exists $v_n\in AC^2_\mathrm{loc}([0,+\infty);\X)$ such that $v_n(0) = x_n$ and $\Vepn(x_n)=\Iepn(v_n)$.

We observe that, for every $T>0$,
\begin{equation}
    \int_0^T|v_n'|(t)^2\,dt\leq
    2\Iepn(v_n)e^{T/\ep_n}=
    2\Vepn(x_n)e^{T/\ep_n},
    \nonumber
\end{equation}
and hence from (\ref{th:AC-Holder}) it follows that
\begin{equation}
    d(v_n(t_2),v_n(t_1))\leq
    \left[2\Vepn(x_n)e^{T/\ep_n}\right]^{1/2}|t_2-t_1|^{1/2}
    \qquad
    \forall (t_1,t_2)\in[0,T]^2.
    \nonumber
\end{equation}

Now we show that
\begin{equation}
    \lim_{n\to +\infty}v_n(\ep_n\sigma)=x_\infty
    \qquad
    \forall\sigma\geq 0.
    \label{claim:vn2x}
\end{equation}
Indeed, from the triangle inequality and the previous estimate with $t_1=0$ and $t_2=T=\ep_n\sigma$, we obtain
\begin{align*}
    d(v_n(\ep_n\sigma),x_\infty)
    &\leq d(v_n(\ep_n\sigma),x_n)+d(x_n,x_\infty)
    \\
    &= d(v_n(\ep_n\sigma),v_n(0))+d(x_n,x_\infty)
    \\
    &\leq \left[2\Vepn(x_n)e^\sigma\right]^{1/2}(\ep_n\sigma)^{1/2}
    +d(x_n,x_\infty),
\end{align*}
so that the claim follows from the boundedness of $\{\Vepn(x_n)\}$ and the convergences $\ep_n\to 0^+$ and $x_n\to x_\infty$.

Finally, from the definition of $\Iepn(v_n)$ we deduce that
\begin{equation}
    \Vepn(x_n)=
    \Iepn(v_n)\geq
    \int_0^{+\infty}F(v_n(t))\frac{e^{-t/\ep_n}}{\ep_n}\,dt=
    \int_0^{+\infty}F(v_n(\ep_n\sigma))e^{-\sigma}\,d\sigma.
    \nonumber
\end{equation}
At this point, Fatou's lemma,~(\ref{claim:vn2x}), and the lower semicontinuity of $F$ yield
\begin{eqnarray*}
    \liminf_{n\to +\infty}\Vepn(x_n) & \geq &
    \liminf_{n\to +\infty}\int_0^{+\infty}F(v_n(\ep_n\sigma))e^{-\sigma}\,d\sigma
    \\
    & \geq &
    \int_0^{+\infty}\liminf_{n\to +\infty}F(v_n(\ep_n\sigma))e^{-\sigma}\,d\sigma
    \\
    & \geq &
    \int_0^{+\infty}F(x_\infty)e^{-\sigma}\,d\sigma
    \\
    & = &
    F(x_\infty),
\end{eqnarray*}
which proves~(\ref{th:liminf-vf}).

\paragraph{\textmd{\textit{Introducing maximal slope curves}}}

In order to prove~(\ref{th:vf-cont}) and~(\ref{th:liminf-slope}), for each $n\in\n$ we consider a maximal slope curve $(w_n,\psi_n,N_n)$ for $F$ with initial condition $w_n(0)=x_n$ (see Theorem~\ref{thmbibl:CMP}). We claim that, whenever the sequence $\{F(x_n)\}$ is bounded,
\begin{equation}
    \lim_{n\to +\infty}w_n(\ep_n\sigma)=x_\infty
    \qquad
    \forall \sigma\geq 0.
    \label{claim:wn2x}
\end{equation}

The proof is analogous to that of~(\ref{claim:vn2x}). Indeed, from the definition of maximal slope curve we obtain
\begin{equation}
    \frac{1}{2}\int_{t_1}^{t_2}|w_n'|(t)^2\,dt\leq
    \psi_n(t_1)-\psi_n(t_2)\leq
    \psi_n(0)=
    F(x_n),
    \nonumber
\end{equation}
and hence, by~(\ref{th:AC-Holder}),
\begin{equation}
    d(w_n(t_2),w_n(t_1))\leq
    \left[2F(x_n)\right]^{1/2}|t_2-t_1|^{1/2}
    \qquad
    \forall (t_1,t_2)\in[0,T]^2.
    \nonumber
\end{equation}

Therefore, applying the previous estimate with $t_1=0$ and $t_2=T=\ep_n\sigma$, and using the triangle inequality, we obtain
\begin{align*}
    d(w_n(\ep_n\sigma),x_\infty)
    &\leq d(w_n(\ep_n\sigma),x_n)+d(x_n,x_\infty)
    \\
    &= d(w_n(\ep_n\sigma),w_n(0))+d(x_n,x_\infty)
    \\
    &\leq \left[2F(x_n)\right]^{1/2}(\ep_n\sigma)^{1/2}
    +d(x_n,x_\infty).
\end{align*}
The claim then follows from the boundedness of $\{F(x_n)\}$ and the convergences $\ep_n\to 0^+$ and $x_n\to x_\infty$.

\paragraph{\textmd{\textit{Convergence of value functions}}}

We show that, under assumptions~(\ref{hp:F-cont}) and~(\ref{hp:bound-slope}), one has
\begin{equation}
    \limsup_{n\to +\infty}\Vepn(x_n)\leq F(x_\infty),
    \label{th:limsup-fv-F}
\end{equation}
which, together with~(\ref{th:liminf-vf}), proves~(\ref{th:vf-cont}).

To this end, we consider the maximal slope curves $(w_n,\psi_n,N_n)$  introduced above. From the theory of maximal slope curves we know that $|w_n'|(t)=|\nabla F|(w_n(t))$ for almost every $t\geq 0$. Therefore, from Lemma~\ref{lemma:CMP1} we obtain that
\begin{align*}
    \Vepn(x_n)
    &\leq
    \int_0^{+\infty}
    \left\{\frac{\ep_n}{2}|w_n'|(\ep_n\sigma)^2+\psi_n(\ep_n\sigma)\right\}
    e^{-\sigma}\,d\sigma
    \\
    &=
    \int_0^{+\infty}
    \left\{\frac{\ep_n}{2}|\nabla F|(w_n(\ep_n\sigma))^2+\psi_n(\ep_n\sigma)\right\}
    e^{-\sigma}\,d\sigma
    \\
    &\leq
    \frac{\ep_n}{2}\cdot\frac{2}{\ep_n}
    \left[F(x_n)-\Vepn(x_n)\right]
    +
    \int_0^{+\infty}\psi_n(\ep_n\sigma)e^{-\sigma}\,d\sigma.
\end{align*}

Due to the second bound in~(\ref{hp:bound-slope}), the first term on the right-hand side tends to~0, and therefore
\begin{equation}
    \limsup_{n\to +\infty}\Vepn(x_n)\leq
    \limsup_{n\to +\infty}
    \int_0^{+\infty}\psi_n(\ep_n\sigma)e^{-\sigma}\,d\sigma.
    \label{est:limsup1}
\end{equation}

We now select an increasing sequence of positive integers $\{n_k\}$ such that
\begin{equation}
    \limsup_{n\to +\infty}
    \int_0^{+\infty}\psi_n(\ep_n\sigma)e^{-\sigma}\,d\sigma
    =
    \lim_{k\to +\infty}
    \int_0^{+\infty}\psi_{n_k}(\ep_{n_k}\sigma)e^{-\sigma}\,d\sigma.
    \label{est:limsup2}
    %\nonumber
\end{equation}

Since the functions $\sigma\mapsto\psi_{n_k}(\ep_{n_k}\sigma)$ are nonincreasing and uniformly bounded, Helly's selection principle implies that, up to a further subsequence, not relabeled, there exists a nonincreasing function $\psi_\infty:[0,+\infty)\to[0,+\infty)$ such that
\begin{equation}
    \lim_{k\to +\infty}\psi_{n_k}(\ep_{n_k}\sigma)=\psi_\infty(\sigma)
    \qquad
    \forall \sigma\geq 0.
\end{equation}
Here the uniform boundedness follows from
\begin{equation}
    0\leq \psi_{n_k}(\ep_{n_k}\sigma)\leq \psi_{n_k}(0)=F(x_{n_k})
    \qquad
    \forall \sigma\geq 0,
    \nonumber
\end{equation}
and from the first bound in~(\ref{hp:bound-slope}). Therefore, by the dominated convergence theorem,
\begin{equation}
    \lim_{k\to +\infty}
    \int_0^{+\infty}\psi_{n_k}(\ep_{n_k}\sigma)e^{-\sigma}\,d\sigma
    =
    \int_0^{+\infty}\psi_\infty(\sigma)e^{-\sigma}\,d\sigma.
    \label{est:limsup3}
\end{equation}

If we show that
\begin{equation}
    \psi_\infty(\sigma)=F(x_\infty)
    \qquad
    \text{for almost every $\sigma\geq 0$,}
    \label{claim:psi-infty}
\end{equation}
then~(\ref{th:limsup-fv-F}) follows from~(\ref{est:limsup1}),~(\ref{est:limsup2}), and~(\ref{est:limsup3}).

In order to prove~(\ref{claim:psi-infty}), we rely on the estimate (\ref{th:vf-vs-slope}) of Lemma~\ref{lemma:CMP1} with $\ep:=\ep_{n_k}$, $x:=x_{n_k}$, and $w:=w_{n_k}$. Recalling the second bound in~(\ref{hp:bound-slope}), we obtain
\begin{equation}
    \sup_{k\in\n}\int_0^{+\infty}|\nabla F|(w_{n_k}(\ep_{n_k}\sigma))^2 e^{-\sigma}\,d\sigma\leq
    \sup_{k\in\n}\frac{2}{\ep_{n_k}}\left[F(x_{n_k})-\Vepnk(x_{n_k})\right]
    <+\infty.
    \nonumber
\end{equation}

As a consequence, Fatou's lemma yields
\begin{multline}
    \qquad
    \int_0^{+\infty}\liminf_{k\to +\infty}|\nabla F|(w_{n_k}(\ep_{n_k}\sigma))^2 e^{-\sigma}\,d\sigma
    \\
    \leq
    \liminf_{k\to +\infty}\int_0^{+\infty}|\nabla F|(w_{n_k}(\ep_{n_k}\sigma))^2 e^{-\sigma}\,d\sigma
    <+\infty,
    \qquad
    \nonumber
\end{multline}
and therefore
\begin{equation}
    \liminf_{k\to +\infty}|\nabla F|(w_{n_k}(\ep_{n_k}\sigma))^2<+\infty
    \qquad
    \text{for almost every $\sigma\geq 0$}.
    \label{en-passant}
\end{equation}

Let $N_\infty$ be a set of measure zero such that the liminf above is finite for every $\sigma\in[0,+\infty)\setminus N_\infty$, and set
\begin{equation}
    S:=N_\infty\cup\bigcup_{k\in\n}\{\sigma\geq 0:\ep_{n_k}\sigma\in N_{n_k}\}.
    \nonumber
\end{equation}
Then $S$ has measure zero. We claim that~(\ref{claim:psi-infty}) holds for every $\sigma\in[0,+\infty)\setminus S$. Indeed, for any such $\sigma$, the liminf in~(\ref{en-passant}) is finite and, by the definition of maximal slope curve,
\begin{equation}
    F(w_{n_k}(\ep_{n_k}\sigma))
    =
    \psi_{n_k}(\ep_{n_k}\sigma)
    \leq
    \psi_{n_k}(0)
    =
    F(x_{n_k}).
    \nonumber
\end{equation}
Therefore, by the first bound in~(\ref{hp:bound-slope}),
\begin{equation}
    \sup_{k\in\n}F(w_{n_k}(\ep_{n_k}\sigma))<+\infty.
    \nonumber
\end{equation}

As a consequence, for any such $\sigma$ there exists a further subsequence $\{n_{k_\ell}\}$, possibly depending on $\sigma$, such that
\begin{equation}
    \sup_{\ell\in\n}F(w_{n_{k_\ell}}(\ep_{n_{k_\ell}}\sigma))<+\infty
    \qquad\text{and}\qquad
    \sup_{\ell\in\n}|\nabla F|(w_{n_{k_\ell}}(\ep_{n_{k_\ell}}\sigma))<+\infty.
    \nonumber
\end{equation}
This is exactly what is needed in order to apply property~(\ref{hp:F-cont}) to the sequence
\begin{equation}
    w_{n_{k_{\ell}}}(\ep_{n_{k_{\ell}}}\sigma)\to x_\infty,
    \nonumber
\end{equation}
whose convergence is due to (\ref{claim:wn2x}). We conclude that
\begin{equation}
    \psi_\infty(\sigma)=
    \lim_{k\to+\infty}\psi_{n_k}(\ep_{n_k}\sigma)=
    \lim_{\ell\to+\infty}\psi_{n_{k_\ell}}(\ep_{n_{k_\ell}}\sigma)=
    \lim_{\ell\to+\infty}F(w_{n_{k_\ell}}(\ep_{n_{k_\ell}}\sigma))=
    F(x_\infty),
    \nonumber
\end{equation}
which completes the proof of (\ref{th:vf-cont}).

\paragraph{\textmd{\textit{Gamma--liminf inequality for slopes}}}

The proof is analogous to that of~(\ref{th:vf-cont}), but with one less level of subsequences. Proceeding as before, we obtain
\begin{equation}
    \liminf_{n\to +\infty}\frac{2}{\ep_n}\left[F(x_n)-\Vepn(x_n)\right]\geq
    \int_0^{+\infty}\liminf_{n\to +\infty}|\nabla F|(w_n(\ep_n\sigma))^2e^{-\sigma}\,d\sigma,
    \nonumber
\end{equation}
so that~(\ref{th:liminf-slope}) follows once we prove that
\begin{equation}
    \liminf_{n\to +\infty}|\nabla F|(w_n(\ep_n\sigma))^2\geq
    |\nabla F|(x_\infty)^2
    \qquad
    \text{for almost every $\sigma\geq 0$.}
    \label{en-passant2}
\end{equation}

Arguing as before, we can show that there exists a set $S$ of measure zero such that, for every $\sigma\geq 0$ with $\sigma\notin S$, there exists a subsequence $\{n_k\}$, possibly depending on $\sigma$, that realizes the liminf of the slopes in~(\ref{en-passant2}) and such that
\begin{equation}
    \sup_{k\in\n}F(w_{n_k}(\ep_{n_k}\sigma))<+\infty
    \qquad\text{and}\qquad
    \sup_{k\in\n}|\nabla F|(w_{n_k}(\ep_{n_k}\sigma))<+\infty.
    \nonumber
\end{equation}

We then conclude by applying property~(\ref{hp:F-cont}) to the sequence
\begin{equation}
    w_{n_k}(\ep_{n_k}\sigma)\to x_\infty,
    \nonumber
\end{equation}
whose convergence is again due to (\ref{claim:wn2x}).
\end{proof}

\begin{rmk}[Equivalent boundedness conditions]
\begin{em}
In Statement~(\ref{statement:deluxe}) of Proposition~\ref{lemma:lim-vf}, the boundedness of $F(x_n)$ in~(\ref{hp:bound-slope}) can equivalently be replaced by the boundedness of $\Vepn(x_n)$. Indeed, under the second condition in~(\ref{hp:bound-slope}), we have
\[
    0\leq F(x_n)-\Vepn(x_n)
    =\frac{\ep_n}{2}\cdot\frac{2[F(x_n)-\Vepn(x_n)]}{\ep_n}
    \longrightarrow 0,
\]
and therefore the two boundedness assumptions are equivalent.
As a consequence, the assumption~(\ref{hp:bound-slope}) and the conclusions~(\ref{th:vf-cont})-(\ref{th:liminf-slope}) have precisely the same structure as the implication~(\ref{hp:slope-energy}) in Theorem~\ref{thm:lgf-gfl}, with the quantity $2(F-\Vep)/\ep$ playing the role of the square of the metric slope of $\Vep$. This provides a further justification for adapting the proof strategy of that theorem to our setting.
\end{em}
\end{rmk}

\begin{rmk}[Gamma-limit of the value function]\label{rmk:Gamma-vf}
\begin{em}
    An immediate consequence of Statement~(\ref{statement:gammaliminf}) of Proposition~\ref{lemma:lim-vf}, together with the trivial estimate $\Vep(x)\leq F(x)$ for every $\ep>0$ and every $x\in\X$, is that $F$ is both the pointwise limit and the $\Gamma$-limit of $\Vep$ as $\ep\to 0^+$. 
    In particular, $\Vep(x)\to F(x)$ from below for every $x\in\X$.
\end{em}    
\end{rmk}

%\clearpage

\subsection{A variant of Fatou's lemma}

The following elementary variant of Fatou's lemma will be useful when a lower bound on the pointwise liminf is available only along subsequences satisfying an additional boundedness condition.

\begin{lemma}\label{lemma:An-Bn}
    Let $(a,b)\subseteq\re$ be an interval. For every $n\in\n$, let $A_n:(a,b)\to[0,+\infty)$ and $B_n:(a,b)\to[0,+\infty)$ be measurable functions, and let $A_\infty:(a,b)\to[0,+\infty)$ be another measurable function.

    Let us assume that 
    \begin{enumerate}
    \renewcommand{\labelenumi}{(\roman{enumi})}
        \item there exist two real numbers $M_A$ and $M_B$ such that, for every $n\in\n$,
        \begin{equation}
            \int_a^b A_n(t)\,dt\leq M_A
            \qquad\quad\text{and}\quad\qquad
            \int_a^b B_n(t)\,dt\leq M_B;
            \label{hp:int-An-Bn}
        \end{equation}

        \item for every $t\in(a,b)$ and every increasing sequence $\{n_k\}$ of positive integers, the following implication holds:
        \begin{equation}
        \fbox{$\begin{array}{c}
            \displaystyle{\sup_{k \in \n}A_{n_k}(t)< +\infty} \\
            \noalign{\vspace{1ex}}
            \displaystyle{\sup_{k \in \n}B_{n_k}(t)< +\infty}
        \end{array}$}
        \quad\Longrightarrow\quad
        \fbox{$\displaystyle\liminf_{k\to +\infty}A_{n_k}(t)\geq A_\infty(t)$}.
        \label{hp:liminf-An}
        \end{equation}
    \end{enumerate}

    Then
    \begin{equation}
        \liminf_{n\to +\infty}\int_a^b A_n(t)\,dt\geq
        \int_a^b A_\infty(t)\,dt.
        \label{th:liminf-int-An}
    \end{equation}
\end{lemma}

\begin{proof}
Let us fix $\delta>0$. By Fatou's lemma and assumption~(\ref{hp:int-An-Bn}), we obtain that
\begin{equation}
    \int_a^b \liminf_{n\to +\infty}[A_n(t)+\delta B_n(t)]\,dt
    \leq
    \liminf_{n\to +\infty}\int_a^b [A_n(t)+\delta B_n(t)]\,dt
    \leq
    M_A+\delta M_B,
    \nonumber
\end{equation}
and hence
\begin{equation}
    \liminf_{n\to +\infty}[A_n(t)+\delta B_n(t)]<+\infty
    \qquad
    \text{for almost every $t\in(a,b)$.}
    \nonumber
\end{equation}

Let us fix such a point $t$. We choose an increasing sequence $\{n_k\}$ of positive integers, possibly depending on $t$, such that
\begin{equation}
    \lim_{k\to +\infty}[A_{n_k}(t)+\delta B_{n_k}(t)]
    =
    \liminf_{n\to +\infty}[A_n(t)+\delta B_n(t)]<+\infty.
    \nonumber
\end{equation}

Since $A_n$ and $B_n$ are nonnegative, this implies that the two sequences $\{A_{n_k}(t)\}$ and $\{B_{n_k}(t)\}$ are bounded from above. Therefore, by~(\ref{hp:liminf-An}),
\begin{equation}
    \liminf_{n\to +\infty}[A_n(t)+\delta B_n(t)]
    =
    \lim_{k\to +\infty}[A_{n_k}(t)+\delta B_{n_k}(t)]
    \geq
    \liminf_{k\to +\infty}A_{n_k}(t)
    \geq
    A_\infty(t).
    \label{est:liminf-An}
    %\nonumber
\end{equation}

On the other hand, by~(\ref{hp:int-An-Bn}),
\begin{equation}
    \int_a^b A_n(t)\,dt+\delta M_B\geq
    \int_a^b [A_n(t)+\delta B_n(t)]\,dt.
\end{equation}

Using Fatou's lemma again, together with~(\ref{est:liminf-An}), we obtain
\begin{align*}
    \delta M_B+\liminf_{n\to +\infty}\int_a^b A_n(t)\,dt
    &\geq
    \liminf_{n\to +\infty}\int_a^b [A_n(t)+\delta B_n(t)]\,dt
    \\
    &\geq
    \int_a^b \liminf_{n\to +\infty}[A_n(t)+\delta B_n(t)]\,dt
    \\
    &\geq
    \int_a^b A_\infty(t)\,dt.
\end{align*}

Since $\delta>0$ is arbitrary, this proves~(\ref{th:liminf-int-An}).
\end{proof}

%\clearpage

\subsection{Limits of minimizers are maximal slope curves}

We are now ready to prove Statement~(\ref{statement:limit}) of Theorem~\ref{thm:main}. The argument follows closely the proof of Theorem~3.1 in~\cite{SlowPM2TV}. To fix the notation, we are given a sequence of positive real numbers $\ep_n\to 0^+$ and, for each $n\in\n$, a minimizer $u_n$ of~(\ref{defn:min-datun}) with $\ep:=\ep_n$. We assume that $u_n$ converges to some limit $u_\infty$ uniformly on every bounded interval $[0,T]$.

We start from the energy identity~(\ref{th:min-CMP}), which now reads as
\begin{multline}
    \qquad
    \frac{1}{2}\int_{t_1}^{t_2}|u_n'|(t)^2\,dt+
    \frac{1}{2}\int_{t_1}^{t_2}\frac{2}{\ep_n}\left[F(u_n(t))-\Vepn(u_n(t))\right]\,dt
    \\
    =
    \Vepn(u_n(t_1))-\Vepn(u_n(t_2))
    \qquad
    \label{eqn:un-CMP}
\end{multline}
for every $0\leq t_1\leq t_2$, and we pass to the limit in both its left-hand side and its right-hand side, which, for brevity, we denote by $\LHS_n(t_1,t_2)$ and $\RHS_n(t_1,t_2)$, respectively.

Our goal is to identify a nonincreasing function $\psi_\infty:[0,+\infty)\to[0,+\infty)$ and a set $N_\infty\subseteq(0,+\infty)$ of measure zero such that
\begin{equation}
    \frac{1}{2}\int_{t_1}^{t_2}|u_\infty'|(t)^2\,dt+
    \frac{1}{2}\int_{t_1}^{t_2}|\nabla F|(u_\infty(t))^2\,dt
    \leq
    \psi_\infty(t_1)-\psi_\infty(t_2)
    \label{th:graal}
\end{equation}
for every $0\leq t_1\leq t_2$, and
\begin{equation}
    \psi_\infty(t)=F(u_\infty(t))
    \qquad
    \forall t\in[0,+\infty)\setminus N_\infty.
    \nonumber
\end{equation}

These two properties are precisely what is needed to show that $(u_\infty,\psi_\infty,N_\infty)$ is a maximal slope curve for $F$.

We denote by $\LHS_\infty(t_1,t_2)$ and $\RHS_\infty(t_1,t_2)$, respectively, the left-hand side and the right-hand side of the expected limit inequality (\ref{th:graal}). We observe that $\LHS_n(t_1,t_2)$ and $\RHS_n(t_1,t_2)$ are uniformly bounded. Indeed, since the value function is nonnegative, we have
\begin{equation}
    \Vepn(u_n(t_1))-\Vepn(u_n(t_2))\leq
    \Vepn(u_n(t_1))\leq
    \Vepn(u_n(0))=
    \Vepn(x_0)\leq
    F(x_0).
    \label{est:LRHS}
\end{equation}

\paragraph{\textmd{\textit{Passing the left-hand side to the limit}}}

We claim that
\begin{equation}
    \liminf_{n\to +\infty}\int_{t_1}^{t_2}|u_n'|(t)^2\,dt\geq
    \int_{t_1}^{t_2}|u_\infty'|(t)^2\,dt
    \label{th:liminf-un'}
\end{equation}
and 
\begin{equation}
    \liminf_{n\to +\infty}\int_{t_1}^{t_2}\frac{2}{\ep_n}[F(u_n(t))-V_{\ep_n}(u_n(t))]\,dt\geq
    \int_{t_1}^{t_2}|\nabla F|(u_\infty(t))^2\,dt,
    \label{th:para-slope}
\end{equation}
and therefore
\begin{equation}
    \liminf_{n\to +\infty}\LHS_n(t_1,t_2)\geq\LHS_\infty(t_1,t_2).
    \nonumber
\end{equation}

The first estimate~(\ref{th:liminf-un'}) is rather standard. Let us consider a subsequence that realizes the liminf. Since this subsequence is bounded in $L^2((t_1,t_2))$, up to passing to a further subsequence it has a weak limit $g$. One can easily check that $g(t)\geq |u_\infty'|(t)$ for almost every $t\in (t_1,t_2)$, and the conclusion follows from the lower semicontinuity of the norm.

In order to show (\ref{th:para-slope}), we apply Lemma~\ref{lemma:An-Bn} in the interval $(t_1,t_2)$ with
\begin{equation*}
    A_n(t):=\frac{2}{\ep_n}[F(u_n(t))-V_{\ep_n}(u_n(t))],
    \qquad\quad
    A_\infty(t):=|\nabla F|(u_\infty(t))^2,
\end{equation*}
and
\begin{equation*}
    B_n(t):=F(u_n(t)).
    \nonumber
\end{equation*}

Let us check the assumptions.
\begin{itemize}
    \item The functions $B_n(t)$ are nonnegative because of the sign condition (\ref{hp:F-sign}) on $F$. The functions $A_n(t)$ are nonnegative because $\Vep(x)\leq F(x)$ for every $x\in\X$. 

    \item The function $t\mapsto|\nabla F|(u_\infty(t))$ is measurable. Indeed, $u_\infty$ is continuous, and the function $x\mapsto|\nabla F|(x)$ is Borel measurable. The latter follows from the identity
    \[
    \left\{x\in\X:|\nabla F|(x)\leq L\strut\right\}
    =
    \bigcup_{n\in\n}
    \left\{x\in\X:|\nabla F|(x)\leq L,\ F(x)\leq n\strut\right\},
    \]
    valid for every $L>0$, since all sets in the union are closed by the assumption~(\ref{hp:F-cont}).
    
    \item The uniform bound on the integrals of $B_n$ follows from the uniform bound (\ref{th:min-bounded}) for minimizers. The uniform bound on the integrals of $A_n$ follows from~(\ref{eqn:un-CMP}) and~(\ref{est:LRHS}).

    \item Assumption (\ref{hp:liminf-An}) follows from Statement~(\ref{statement:deluxe}) in Proposition~\ref{lemma:lim-vf} applied to the sequences $\ep_{n_k}\to 0^+$ and $u_{n_k}(t)\to u_\infty(t)$.
\end{itemize}

\paragraph{\textmd{\textit{Identifying the limit of the right-hand side}}}

We observe that, for every $n\in\n$, the function $t\mapsto\Vepn(u_n(t))$ is nonincreasing and nonnegative. Hence, by Helly's selection principle, there exist an increasing sequence of positive integers $\{n_k\}$ and a nonincreasing function $\psi_\infty:[0,+\infty)\to[0,+\infty)$ such that
\begin{equation}
    \lim_{k\to+\infty}\Vepnk(u_{n_k}(t))=
    \psi_\infty(t)
    \qquad  
    \forall t\geq 0.
    \label{defn:psi}
\end{equation}

Even though $\psi_\infty$ is obtained only along a subsequence, this is enough to conclude~(\ref{th:graal}), because
\begin{align*}
    \LHS_\infty(t_1,t_2)
    &\leq
    \liminf_{n\to +\infty}\LHS_n(t_1,t_2)
    =
    \liminf_{n\to +\infty}\RHS_n(t_1,t_2)
    \\
    &\leq
    \liminf_{k\to +\infty}\RHS_{n_k}(t_1,t_2)
    =
    \RHS_\infty(t_1,t_2).
\end{align*}

\paragraph{\textmd{\textit{Characterization of the limit of the right-hand side}}}

We claim that
\begin{equation}
    \psi_\infty(t)=F(u_\infty(t))
    \qquad
    \text{for almost every $t\geq 0$, including $t=0$}.
    \label{eq:char-psi}
\end{equation}

To this end, we consider the sequence $\{n_k\}$ involved in the definition~(\ref{defn:psi}) of $\psi_\infty$, and we fix $T>0$. From~(\ref{eqn:un-CMP}),~(\ref{est:LRHS}), and~(\ref{th:min-bounded}) we obtain
\begin{eqnarray*}
    \lefteqn{\hspace{-3em}
    \int_0^T\liminf_{k\to+\infty}
    \left\{\frac{2}{\ep_{n_k}}\left[F(\unk(t))-\Vepnk(\unk(t))\right]+
    F(\unk(t))\right\}dt}
    \\[0.5ex]
    & \leq &
    \liminf_{k\to+\infty}\int_0^T\left\{\frac{2}{\ep_{n_k}}\left[F(\unk(t))-\Vepnk(\unk(t))\right]+
    F(\unk(t))\right\}dt
    \\[1ex]
    & \leq &
    (T+3)F(x_0),
\end{eqnarray*}
and hence
\begin{equation}
    \liminf_{k\to+\infty}
    \left[\frac{2}{\ep_{n_k}}\left\{F(\unk(t))-\Vepnk(\unk(t))\right\}+
    F(\unk(t))\right]<+\infty
    \label{est:limin-nk}
\end{equation}
for almost every $t\in[0,T]$. Since $T>0$ is arbitrary,~(\ref{est:limin-nk}) holds for almost every $t\geq 0$.

We claim that~(\ref{eq:char-psi}) holds for every such $t$. Indeed, for every such $t$, there exists a further subsequence $\{n_{k_\ell}\}$, possibly depending on $t$, such that
\begin{equation}
    \sup_{\ell\in\n}F(u_{n_{k_\ell}}(t))<+\infty
    \qquad\text{and}\qquad
    \sup_{\ell\in\n}\frac{2}{\ep_{n_{k_\ell}}}
    \left[F(u_{n_{k_{\ell}}}(t))-V_{\ep_{n_{k_{\ell}}}}(u_{n_{k_{\ell}}}(t))\right]<+\infty.
    \nonumber
\end{equation}

Applying Statement~(\ref{statement:deluxe}) of Proposition~\ref{lemma:lim-vf} to the sequences
\begin{equation}
    \ep_{n_{k_\ell}}\to 0
    \qquad\quad\text{and}\quad\qquad
    u_{n_{k_\ell}}(t)\to u_\infty(t),
    \nonumber
\end{equation}
we conclude that
\begin{equation}
    \psi_\infty(t)=
    \lim_{k\to+\infty}\Vepnk(u_{n_k}(t))=
    \lim_{\ell\to+\infty}V_{\ep_{n_{k_{\ell}}}}(u_{n_{k_{\ell}}}(t))=
    F(u_\infty(t)).
    \nonumber
\end{equation}

It remains to show that the equality in~(\ref{eq:char-psi}) holds for $t=0$, but this follows from Remark~\ref{rmk:Gamma-vf} because
\[
    \psi_\infty(0)
    =
    \lim_{k\to+\infty}\Vepnk(u_{n_k}(0))
    =
    \lim_{k\to+\infty}\Vepnk(x_0)
    =
    F(x_0).
\]

This completes the proof.
\qed

\setcounter{equation}{0}
\section{Open problems}\label{sec:open}

We conclude the paper by mentioning two questions that naturally arise from our analysis and might deserve further investigation.

\begin{open}[The role of completeness] \label{open:compl}

Is the completeness of the metric space $\X$ really necessary in Theorem~\ref{thm:main}?

\end{open}

The classical existence theory for maximal slope curves does not require $\X$ to be complete, since the compactness needed in the construction is already provided by the local coercivity assumption~(\ref{hp:F-loc-coerc}). In the present paper, completeness is used only in Step~3 of the proof of the compactness criterion for $AC^2$ trajectories (Proposition~\ref{prop:cpt}). Without completeness, Proposition~\ref{prop:cpt} is no longer true. However, this does not exclude the possibility that statements~(\ref{statement:existence}) and~(\ref{statement:compactness}) of Theorem~\ref{thm:main} remain valid even without completeness.

\begin{open}[Metric slope of the value function]

Under which assumptions is it true that
\begin{equation}
    |\nabla\Vep|(x)^2=
    \frac{2}{\ep}\left[F(x)-\Vep(x)\right]
    \label{open:slope}
\end{equation}
for every $x\in\X$ such that $\Vep(x) < \infty$?

\end{open}

When~(\ref{open:slope}) holds, every minimizer $\uep$ is itself a maximal slope curve for the value function $\Vep$. Indeed, in that case, the energy identity~(\ref{th:min-CMP}) coincides with the defining identity of a maximal slope curve. Moreover, together with the convergence properties established in Proposition~\ref{lemma:lim-vf}, this would allow us to apply Theorem~\ref{thm:lgf-gfl} directly, rather than adapting its proof.

As a consequence of~\cite[Theorem~7.1]{rossi_savare_segatti_stefanelli_2019}, identity~(\ref{open:slope}) holds when~$\X$ is a geodesic space, the conditions~(\ref{hp:F-sign})-(\ref{hp:F-lsc})-(\ref{hp:F-loc-coerc}) hold, and, additionally,~$F$ is geodesically semiconvex and continuous. Since the metric slope is a local quantity, it would be natural, for example, to attempt relaxing the semiconvexity assumption.

\subsubsection*{\centering Acknowledgments}
The first author is a member of the Italian {\selectlanguage{italian} ``Gruppo Nazionale per l'Analisi Matematica, la Probabilità e le loro Applicazioni'' (GNAMPA) of the ``Istituto Nazionale di Alta Matematica'' (INdAM)} and acknowledges the MIUR Excellence Department Project awarded to the Department of Mathematics of the University of Pisa (CUP I57G22000700001). 

The second author acknowledges the support of the French ``Agence nationale de la recherche'' (ANR) through the PEPR PDE-AI project (ANR-23-PEIA-0004).

\end{document}